\documentclass[12pt,reqno]{amsart}

\usepackage{amsfonts,color,amsthm,amsmath,amssymb}

\usepackage{color}
\def\Box{\vcenter{\vbox{\hrule\hbox{\vrule
     \vbox to 8.8pt{\hbox to 10pt{}\vfill}\vrule}\hrule}}}

\newcommand{\tr}{\text{Tr}}
\def\qed{{\hfill$\square$}}

\def\Z{{\mathbb Z}}

\def\C{{\mathbb C}}
\def\F{{\mathbb F}}

\def\Tr{{\mathrm{Tr}}}
\def\Norm{{\mathrm{Norm}}}
\def\Cay{{\mathrm{Cay}}}

\newcommand{\ga}{{\gamma}}

\newtheorem{thm}{Theorem}

\newtheorem{lem}[thm]{Lemma}
\newtheorem{remark}[thm]{Remark}

\numberwithin{equation}{section}

\begin{document}
\title{Three-class association schemes from cyclotomy}
\author{Tao Feng and Koji Momihara  }
\thanks{T. Feng is with Department of Mathematics, Zhejiang University, Hangzhou 310027, Zhejiang, China
(e-mail: tfeng@zju.edu.cn). The work of T. Feng was supported in part by the Fundamental Research Funds for the Central
Universities of China, Zhejiang Provincial Natural Science Foundation (LQ12A01019), National Natural Science
Foundation of China (11201418). }
\thanks{K. Momihara is with Faculty of Education, Kumamoto University,
2-40-1 Kurokami, Kumamoto 860-8555, Japan (e-mail:
momihara@educ.kumamoto-u.ac.jp). The work of
K. Momihara was supported by JSPS under Grant-in-Aid for Research
Activity Start-up 23840032.}
\maketitle

\begin{abstract} We give three constructions of three-class association schemes as fusion schemes of the cyclotomic scheme,
 two of which are primitive.
\end{abstract}

\section{Introduction}
Association schemes form a central part of algebraic combinatorics, and plays important roles in several branches of mathematics, such as coding theory and graph theory. Two-class symmetric association schemes are equivalent to strongly regular graphs, and are extensively studied. The natural graph theoretical extension of strongly regular graph is distance regular graph, whose distance relations form an association scheme. Distance regular graphs have attracted considerable attention, and important progress has been achieved on this topic. We refer the reader to the book \cite{bcn} and the undergoing survey \cite{vkt}. There are not so many papers about three-class association schemes, see the survey \cite{3c} by van Dam and the references therein.
It is the purpose of this note to present new constructions of
primitive three-class association schemes using cyclotomy in finite fields.
As consequences, we obtain three new infinite families of three-class association schemes, two of which are primitive.

Quite recently, there have been several constructions of strongly regular graphs with new parameters and skew Hadamard difference sets from cyclotomy, the latter giving rise to two-class nonsymmetric association schemes, see \cite{FX111,FMX11,GXY11,Momi} for strongly regular graphs and \cite{CF2,FX113,FMX11,M2} for skew Hadamard difference sets.
In \cite{sw}, the authors discussed the problem when a Cayley graph on a finite field with a single cyclotomic class as its connection set can form a  strongly regular graph. Such a strongly regular graph is called {\it cyclotomic}.
They raised the following conjecture: if the Cayley graph on the finite field $\F_q$ of order $q=p^f$ with a multiplicative subgroup $C$ of index $M$ of $\F_q$ as its connection set is cyclotomic strongly regular, then either of the following holds:
\begin{enumerate}
\item[(1)] (subfield case) $C$ is the multiplicative group of a subfield of $\F_{q}$,
\item[(2)] (semi-primitive case) $-1\in \langle p\rangle\le \Z_M^\ast$,
\item[(3)] (exceptional case) it is either of eleven sporadic examples of
cyclotomic strongly regular graphs (see \cite[Table~1]{sw}).
\end{enumerate}
This conjecture is still open but the authors gave a proof in a partial case assuming
the generalized Riemann hypothesis. On the other hand, in \cite{FX111,FMX11,GXY11,M2},
several of these sporadic examples have been generalized into infinite families by taking a union of several cyclotomic classes and doing detailed computations using Gauss sums. For other constructions of strongly regular graphs from cyclotomy, we refer the reader to the references in \cite{FMX11}.

Also, recently, skew Hadamard difference sets are currently under intensive study. There was a major
conjecture in this area:  Up to equivalence the Paley (quadratic residue) difference sets
are the only skew Hadamard difference sets in abelian groups.
This conjecture turned out to be false:
Ding
and Yuan \cite{DY06} gave two counterexamples of this conjecture in finite fields with characteristic three. Furthermore,
Muzychuk \cite{M} constructed infinitely many inequivalent skew Hadamard
difference sets in elementary abelian groups of order $q^3$.
Recently, in \cite{CF2,FX113,FMX11,M2}, the authors constructed further counterexamples of this conjecture by taking suitably a union of  cyclotomic classes. See the introduction of \cite{FX113} (or \cite{CF2})
for a short survey on skew Hadamard difference sets.

Thus, a lot of strongly regular graphs and skew Hadamard difference sets have been obtained from cyclotomy. 
Therefore, we can say that  the cyclotomy  is a quite powerful tool to construct two-class association schemes.
In this note, we shall try to construct three-class association schemes from cyclotomy involving computations
of Gauss sums based on the Hasse-Davenport theorem.

This note is organized as follows:
In Section~\ref{sec:pre}, we review about association schemes and characters of finite fields. In Section~\ref{sec:par}, we
introduce a partition of $\Z_M$, and compute some group ring elements in $\Z[\Z_M]$ based on the results in \cite{ADJP}. 
In Section~\ref{sec:con}, we give three constructions of three-class
association schemes in finite fields with characteristic $2$ as fusion schemes of the cyclotomic schemes.
The parameters of association schemes obtained in Section~\ref{sec:con} are listed in the appendix. We shall use the standard notations on group rings as can be found in the book \cite{BJL}.

\section{Preliminaries}\label{sec:pre}
Let $X$ be a nonempty finite set, and a set of symmetric relations $R_0,R_1,\cdots, R_d$ be a partition of $X\times X$ such that $R_0=\{(x,x)|x\in X\}$. Denote by $A_i$ the adjacency matrix of $R_i$ for each $i$, whose $(x,y)$-th entry is $1$ if $(x,y)\in R_i$ and $0$ otherwise. We call $(X,\{R_i\}_{i=0}^d)$ a {\it $d$-class association scheme} if there exist numbers $p_{i,j}^k$ such that
\[
A_iA_j=\sum_{k=0}^dp_{i,j}^kA_k.
\]
These numbers are called the intersection numbers of the scheme. The $\C$-linear span of $A_0,A_1,\cdots,A_d$ forms a semisimple  algebra of dimension $d+1$, called the {\it Bose-Mesner algebra} of the scheme. With respect to the basis $A_0,A_1,\cdots,A_d$, the matrix of the multiplication by $A_i$ is denoted by $B_i$, namely
\[
A_i(A_0,A_1,\cdots,A_d)=(A_0,A_1,\cdots,A_d)B_i,\;0\leq i\leq d.
\]
Since each $A_i$ is symmetric, this algebra is commutative. There exists a set of minimal idempotents $E_0,E_1,\cdots,E_d$  which also forms a basis of the algebra. The $(d+1)\times (d+1)$ matrix $P$ such that
\[
(A_0,A_1,\cdots,A_d)=(E_0,E_1,\cdots,E_d)P
\]
is called the {\it first eigenmatrix} of the scheme. Dually, the $(d+1)\times (d+1)$ matrix $Q$ such that
\[
(E_0,E_1,\cdots,E_d)=\frac{1}{|X|}(A_0,A_1,\cdots,A_d)Q
\]
is called the {\it second eigenmatrix} of the scheme. We clearly have $PQ=|X|I$.\\

We call an association scheme $(X,\{R_i\}_{i=0}^d)$ a {\it translation association scheme} or a {\it Schur ring} if $X$ is a  (additively
written) finite abelian  group and there exists a partition $S_0=\{0\},S_1,\cdots,S_d$  of $X$ such that
\[
R_i=\{(x,x+y)|\,x\in X, y \in S_i\}.
\]
For brevity, we will just say that $(X,\{S_i\}_{i=0}^d)$ is an association scheme. \\

Assume that $(X,\{S_i\}_{i=0}^d)$ is a translation association scheme. There is an equivalence relation defined on the character group $\hat{X}$ of $X$ as follows: $\chi\sim\chi'$ if and only if $\chi(S_i)=\chi'(S_i)$ for each $0\leq i\leq d$. Here $\chi(S)=\sum_{g\in S}\chi(g)$, for any $\chi\in\hat{X}$, and $S\subseteq X$. Denote by $D_0, D_1,\cdots,D_d$ the equivalence classes, with $D_0$ consisting of only the principal character. Then $(\hat{X},\{D_i\}_{i=0}^d)$ forms a translation association scheme, called the {\it dual} of $(X,\{S_i\}_{i=0}^d)$. The first eigenmatrix of the dual scheme is equal to the second eigenmatrix of the original scheme. Please refer to \cite{bi} and \cite{bcn} for more details.\\

A classical example of translation schemes is the cyclotomic scheme which we describe now. Let $p$ be a prime and $q=p^f \,(f\geq 1)$ be a prime power, $M|q-1$, and $\ga$ be a primitive element of the finite field $F=\F_q$. Define the multiplicative subgroup $C_0^{(M,F)}=\langle\ga^M\rangle$. Its cosets $C_i^{(M,F)}=\ga^iC_0^{(M,F)}$, $0\leq i\leq M-1$, are called the {\it cyclotomic classes} of order $M$ of $F$. Together with $\{0\}$, they form an $M$-class association scheme, which is called the {\it cyclotomic scheme}. To describe its first eigenmatrix, we define the {\it Gauss periods}
\[
\eta_a=\sum_{x\in C_a^{(M,F)}}\psi(x),\;0\leq a\leq M-1,
\]
where $\psi$ is the canonical additive character of $F$ defined by $\psi(x)=e^{\frac{2\pi i}{p}\tr(x)}$, $x\in F$. The first eigenmatrix $P$ of the scheme is
\begin{equation}\label{eigenmatrix}
P=\left(\begin{array}{cccccc}
1& \frac{q-1}{M}&\frac{q-1}{M}&\frac{q-1}{M}&\cdots &\frac{q-1}{M}\\
1&\eta_{0}  &\eta_1   &\eta_2 & \cdots     &\eta_{M-1}   \\
1&\eta_{1}  &\eta_{2} & \eta_3& \cdots  &\eta_{0} \\
\vdots & & & & \\
1&\eta_{M-1} &\eta_{0} &\eta_1& \cdots  &\eta_{M-2}\\
\end{array}\right).
\end{equation}

For each multiplicative character $\chi$ of $F_q^\ast$, the multiplicative group of $F$, we define the {\it Gauss sum}
\[
G_{F }(\chi)=\sum_{x\in F^\ast}\psi(x)\chi(x).
\]
The following relation will be repeatedly used in this paper (cf. \cite[P.~195]{LN97}):
\[
\psi(x)=\frac{1}{q-1}\sum_{\chi\in\widehat{F^\ast}}G_F(\chi)\chi^{-1}(x),\;\forall\, x\in F^\ast.
\]
Then, the Gauss period can be expressed as a linear combination of Gauss sums as follows:
\begin{align*}
\eta_i&=\psi( \ga^i C_0^{(M,F)})\\
&=\frac{1}{q-1}\sum_{\chi\in\widehat{F^\ast}}G_{F }(\chi)\chi^{-1}(\ga^i)\sum_{x\in C_0^{(M,F)}}\chi^{-1}(x)\\
&=\frac{1}{M}
\sum_{i=0}^{M-1}G_F(\phi^{-i})\phi(\ga^{i} ),
\end{align*}
where $\phi$ is a multiplicative character of order $M$ of $F^\ast$.\\

In this note, we are interested in the fusion schemes of the cyclotomic scheme, namely schemes whose relations are unions of the relations in the cyclotomic scheme. We shall need the following well-known criterion due to Bannai \cite{BannaiSub} and Muzychuk \cite{Muzthesis}, called the {\it Bannai-Muzychuk criterion}: {\it Let $P$ be the first eigenmatrix of an association scheme $(X, \{R_i\}_{0\leq i\leq d})$, and $\Lambda_0:=\{0\}, \Lambda_1,\ldots ,\Lambda_{d'}$ be a partition of $\{0,1,\ldots ,d\}$. Then $(X, \{R_{\Lambda_i}\}_{0\leq i\leq d'})$ forms an association scheme if and only if there exists a partition $\{\Delta_i\}_{0\leq i\leq d'}$ of $\{0,1,2,\ldots ,d\}$ with $\Delta_0=\{0\}$ such that each $(\Delta_i, \Lambda_j)$-block of $P$ has a constant row sum. Moreover, the constant row sum of the $(\Delta_i, \Lambda_j)$-block is the $(i,j)$-th entry of the first eigenmatrix of the fusion scheme.} \\

We close this section by recording the well-known {\it Hasse-Davenport theorem} on Gauss sums. 
\begin{thm}\label{thm:lift}(\cite[Theorem~11.5.2]{BEW97})
Let $\chi$ be a nonprincipal multiplicative character of $\F_q=\F_{p^f}$ and
let $\chi'$ be the lifted character of $\chi$ to the extension field $\F_{q'}=\F_{p^{fs}}$, i.e., $\chi'(\alpha):=\chi(\Norm_{\F_{q'}/\F_q}(\alpha))$ for any $\alpha\in \F_{q'}^*$.
Then, it holds that
\[
G_{\F_{q'}}(\chi')=(-1)^{s-1}(G_{\F_q}(\chi))^s.
\]
\end{thm}

\section{A partition of $\Z_M$, $M=\frac{2^{3s}-1}{2^s-1}$}\label{sec:par}
Let $s$ be a positive integer, and set $M:=\frac{2^{3s}-1}{2^s-1}$. Denote by $F:=\F_{2^{3s}}$,\; $E:=\F_{2^{s}}$ the finite field with 
$2^{3s}$ and $2^s$ elements respectively.  Let
\begin{equation}\label{inv-trace0}
D:=\{ u\in F^*:\;\tr_{F/E}(u^{-1})=0\}.
\end{equation}
This set $D$ is $E^\ast$-invariant, namely $Dg=\{dg:\, d\in D\}=D$ for any $g\in E^\ast$. Therefore $\psi(\omega^a D)$
 depends only on $a\pmod{M}$. First, we show that
$\psi(\omega^a D)$, $0\le a\le M-1$, take exactly three values.
Since $D$ is a union of $E^*$ cosets, we have
\begin{align*}
\psi(\omega^aD)&=\frac{1}{2^s-1}\sum_{u\in D}\sum_{x\in E^*}\psi(x\omega^a u)\\
&=\#\{u:\,u\in D,\,\tr_{F/E}(\omega^a u)=0\}-\frac{1}{2^s-1}\#\{u:\,u\in D,\,\tr_{F/E}(\omega^au)\ne 0\}\\
&=-(2^s+1)+\frac{2^s}{2^s-1}\#\{u:\,u\in D,\,\tr_{F/E}(\omega^au)=0\}.
\end{align*}
It is clear that $u^{\frac{2^{3s}-1}{2^s-1}} \Tr_{F/E}(u^{-1})=\Tr_{F/E}(u^{1+2^s})$, so
\[
D=\{  u\in F^*:\;\tr_{F/E}(u^{1+2^s})=0\}.
\]
Since $Q(x)=\tr_{L/F}(x^{1+2^s})$  is a
nondegenerate quadratic form, the corresponding quadric $\mathcal{Q}$ in $PG(2,2^s)$ intersects $2^s+1$ lines in $1$ point, and $2^{2s-1}+2^{s-1}$ lines in $2$ points \cite{Game}. According to \cite[p.~328]{ADJP}, the tangent lines are given by
\[
L_a=\{[x]:\,x\in F^*,\,\tr_{F/E}(ax)=0\}
\]
with $\tr_{F/E}(a)=0$, $a\ne 0$. Here we  use $[x]$ for the projective point corresponding to the $1$-dimensional subspace spanned by $x$ for each $x\in F^*$.

Therefore, the set
 \[
 S_a:=\{u:\,\tr_{F/E}(u^{1+2^s})=0,\,\tr_{F/E}(\omega^au)=0\}
 \]
has size $0$, $2^s-1$ or $2(2^s-1)$, depending on whether $L_{w^a}$ is a passant line, a tangent line or a secant line. Denote by $T_1$ (resp. $T_2$) those $a$ in $\Z_M$ such that $S_a$ has size $2^s-1$ (resp. $2(2^s-1)$). Then $|T_1|=2^s+1$, $|T_2|=2^{2s-1}+2^{s-1}$. Denote by $T_3$ the remaining elements of $\Z_M$ other than $T_1$ and $T_2$.
We have $|T_3|=2^{2s-1}-2^{s-1}$. To sum up, we have the following result.
\begin{lem}\label{lem:base} The sums $\psi(\omega^a D)$, $0\le a\le M-1$, take exactly three
values, which are
\[
\psi(\omega^a D)=\begin{cases}-1&\;\textup{ if } a\in T_1,\,\\
2^s-1&\;\textup{ if } a\in T_2,\\-2^s-1&\;\textup{ if } a\in T_3.\end{cases}
\]
\end{lem}
The sets $T_1$, $T_2$ and $T_3$ form a partition of $\Z_M$. 
We now prove the following lemma, which is essential for our construction.
\begin{lem}\label{gr_t1} With the above notations, we have
\begin{align}
(T_2-T_3)T_1^{(-1)}&=2^sT_1,\label{ge_t1}\\
(T_2-T_3)T_2^{(-1)}&=2^{2s-1}+2^{s-1}(\Z_M-T_1),\label{ge_t2}\\
(T_2-T_3)T_3^{(-1)}&=-2^{2s-1}+2^{s-1}(\Z_M-T_1).\label{ge_t3}
\end{align}
\end{lem}
\noindent{\bf Proof:}  First we recall from \cite[p.~327]{ADJP} that
\[
T_1=\{i\in \Z_M:\,\tr_{F/E}(w^i)=0\}
\]
by examining the tangent lines of the quadric $\mathcal{Q}$. It is clear that $T_1$ is the classical Singer difference set in $\Z_M$,  so it holds in the group ring $\Z[\Z_M]$ that (see \cite{BJL})
\begin{equation}\label{t1t1}
T_1T_1^{(-1)}=2^s+\Z_M.
\end{equation}
Moreover, we observe that  $\{2i:\,i\in T_1\}$ is equal to $T_1$ by the definition of $T_1$.

We first show that
\[
T_1^2=T_1+2T_2.
\]
For any $i,j\in T_1$, the line
\[
L_{\omega^{i+j}}:=\{[x]:\, x\in F^*,\,\tr_{F/E}(\omega^{i+j}x)=0\}
\]
intersects the quadric  $\mathcal{Q}$ at the points $[\omega^{-i}]$ and $[\omega^{-j}]$, since the set of equations
\[
\tr_{F/E}(X^{-1})=0, \quad \tr_{F/E}(\omega^{i+j}X)=0
\]
has the solutions $X=\omega^{-i}$, $X=\omega^{-j}$. It follows that if $i,j$ are two distinct elements in $\Z_M$, then these two points are distinct and $L_{w^{i+j}}$ is a secant line. Together with  $\{2i:\,i\in T_1\}=T_1$,  we conclude that in $T_1^2$, each element of $T_1$ has coefficient $1$,  each element of $T_3$ has coefficient $0$, and each element of $T_2$ has even coefficient. Now write $d_i$ for the coefficients of $i$ in $T_1^2$ for each $i\in\Z_M$. By direct computation, we have $(T_1T_1^{(-1)})^2=2^{2s}+(2^{2s}+3\cdot 2^s+1)\Z_M$. Examining the coefficient of the identity on both sides of the equation, we get
\[
|T_1|\cdot 1^2+\sum_{i\in T_2}d_i^2=2^{2s+1}+3\cdot2^s+1,
\]
which yields $\sum_{i\in T_2}d_i^2=4|T_2|$. Also, we have $\sum_{i\in T_2}d_i=|T_1|^2-|T_1|=2|T_2|$. Since $d_i/2$  is a nonnegative integer for each $i\in T_2$, and
\[
\sum_{i\in T_2}\left(\frac{d_i}{2}\right)^2=|T_2|,\quad \sum_{i\in T_2}\left(\frac{d_i}{2}\right)=|T_2|,
\]
we immediately get $d_i=2$ for any $i\in T_2$.

We have $T_1+2T_2=\Z_M+(T_2-T_3)$, $T_1^{(-1)}\Z_M=(2^s+1)\Z_M$. Multiplying both sides of $T_1^2=T_1+2T_2$ with $T_1^{(-1)}$, we get
\[
T_1\cdot(2^s+\Z_M)=(2^s+1)\Z_M+(T_2-T_3)T_1^{(-1)}.
\]
The Eqn. \eqref{ge_t1} then follows.

Since $T_1+T_2+T_3=\Z_M$, Eqn. \eqref{ge_t1} yields that $(T_2-T_3)(\Z_M-T_2-T_3)^{(-1)}=2^s T_1$. On the other hand,
 $(T_2-T_3)(T_2-T_3)^{(-1)}=2^{2s}$ by \cite[p.~328]{ADJP}.
Combining these equations, we get Eqn. \eqref{ge_t2} and Eqn. \eqref{ge_t3}.
\qed

\begin{remark} We deduce from Eqn. \eqref{ge_t1}, Eqn. \eqref{t1t1} and $T_1+T_2+T_3=\Z_M$ that
\begin{align}
T_1T_2^{(-1)}&=2^{s-1}T_1^{(-1)}+2^{s-1}\Z_M-2^{s-1},\label{t1t2}\\
T_1T_3^{(-1)}&=-2^{s-1}T_1^{(-1)}+2^{s-1}\Z_M-2^{s-1}.\label{t1t3}
\end{align}
The following equations then follow from direct computations:
\begin{align}
T_1^2T_1^{(-1)}&=2^sT_1+(2^s+1)\Z_M, \label{ttt_inv}\\
T_1^2T_2^{(-1)}& =2^{2s-1}+(2^{s-1}+2^{2s-1})\Z_M-2^{s-1}T_1,\label{ttt_inv2}\\
T_1^2T_3^{(-1)}&=-2^{2s-1}+2^{2s-1}\Z_M-2^{s-1}T_1\label{ttt_inv3}.
\end{align}
\end{remark}
\section{Three-Class Association Schemes in $\F_{2^{3s}}$ and Their Extensions to $\F_{2^{6s}}$ and $\F_{2^{9s}}$}
\label{sec:con}
We fix the following notations throughout this section:
Let $s$ be a positive integer, $M=\frac{2^{3s}-1}{2^s-1}$, and let
$T_1$, $T_2$, $T_3$ be as introduced in the previous section. We define
\[
H:=\F_{2^{9s}},\; G:=\F_{2^{6s}},\; F:=\F_{2^{3s}},\; E:=\F_{2^{s}}.
\]
Let $C_i^{(M,F)}$, $C_i^{(M,G)}$, $C_i^{(M,H)}$, $0\leq i\leq M-1$, be the cyclotomic classes of order $M$ in $F$, $G$, $H$ respectively.
Clearly $C_0^{(M,F)}$ is equal to $E^*$, the multiplicative group of $E$.
Let $\psi$, $\psi'$ $\psi''$ be the canonical additive character of $H$,  $G$ and $F$ respectively.   Also, write
 $\eta_a$, $\eta_a'$, $\eta_a''$, $0\leq a\leq M-1$ for their Gauss periods respectively. Fix a primitive element $\beta$ of $H$ and a primitive element $\gamma$ of $G$ such that $\Norm_{H/F}(\beta)=\Norm_{G/F}(\gamma)$, where $\Norm_{H/F}$ and $\Norm_{G/F}$ is the norm from $H$ to $F$ and from $G$ to $F$ respectively. Write
\[
\omega:=\Norm_{H/F}(\beta)=\Norm_{G/F}(\gamma),
\]
which is a primitive element of $F$.

\subsection{Imprimitive Association Schemes in $\F_{2^{3s}}$}
In this section, we construct an imprimitive three-class association scheme in $\F_{2^{3s}}$. 
Now we prove the following theorem.
\begin{thm}\label{first-ass}
Take the following partition of $F$:
\[
R_0=\{0\},\;R_1=\bigcup_{i\in T_1} C_i^{(M,F)},\; R_2=\bigcup_{i\in T_2} C_i^{(M,F)},\;R_3=\bigcup_{i\in T_3} C_i^{(M,F)}.
\]
Then, $(F,\{R_i\}_{i=0}^3)$ is a three-class association scheme, whose parameters are listed in the appendix.
\end {thm}

\noindent{\bf Proof of Theorem \ref{first-ass}:}
As before, it is easily verified that $\psi(\omega^a C_0^{(M,F)})=2^s-1$ or $-1$ according to $\tr_{F/E}(\omega^a)=0$ or not, i.e.,  $a\in T_1$ or
$a\not \in T_1$. Now,
we compute that
\begin{align*}
\psi(\omega^a R_k)&=\sum_{i\in T_k}\psi(\omega^{a+i}C_0^{(M,F)})\\
&=(2^s-1)|T_1\cap (a+T_k)|-|(\Z_M\setminus T_1)\cap (a+T_k)|\\
&=2^s |T_1\cap (a+T_k)|-|T_k|.
\end{align*}
The term $|T_1\cap (a+T_k)|$ is the coefficient of $a$ in the group ring element $T_1T_k^{(-1)}$. We have computed
 $T_1T_k^{(-1)}$, $1\leq k\leq 3$, in Eqn. \eqref{t1t1}-\eqref{t1t3}.  
For each $k=1,2,3$, the sum $\psi(\omega^a R_k)$ is now computed directly and listed in Table \ref{tab_1}.
By the Bannai-Muzychuk criterion, $(F,\{R_i\}_{i=0}^3)$ is a three-class association scheme.\qed
\vspace{0.3cm}

\begin{table}
\begin{center}
\caption{\label{tab_1}The  values of $\psi(\omega^a R_k)$'s}
\begin{tabular}{|c|c|c|c|c|}\hline
   & $R_0$& $R_1$& $R_2$ &$R_3$\\\hline
$\omega^a=0$ & $1$ & $2^{2s}-1$ & $2^{s-1}(2^{2s}-1)$ & $2^{s-1}(2^s-1)^2$\\\hline
 $a=0$& $1$ & $2^{2s}-1$ &  $-2^{s-1}(2^s+1)$ & $-2^{s-1}(2^s-1)$\\\hline
 $a\in -T_1$ & $1$& $-1$ & $2^{s-1}(2^s-1)$ &   $-2^{s-1}(2^s-1)$\\\hline
 $a\not \in -T_1\cup \{0\}$& $1$ & $-1$ &$-2^{s-1}$ & $2^{s-1}$\\\hline
\end{tabular}
\end{center}
\end{table}

\begin{remark}
By the proof of Theorem~\ref{first-ass}, the dual scheme of the association scheme in Theorem~\ref{first-ass} is given by
\[
D_0=\{0\},\;D_1=C_0^{(M,F)},\; D_2=\cup_{i\in -T_1} C_i^{(M,F)},\;D_3=\cup_{i\in \Z_M\setminus (-T_1\cup \{0\})} C_i^{(M,F)}.
\]
This scheme is imprimitive since $D_0\cup D_1=E$. Their character values are listed in Table \ref{tab_1d}, which we shall need later.
Observe that $D_2=D$.
\begin{table}
\begin{center}
\caption{\label{tab_1d}The  values of $\psi(\omega^a D_k)$'s}
\begin{tabular}{|c|c|c|c|c|}\hline
   & $D_0$& $D_1$& $D_2$ &$D_3$\\\hline
$\omega^a=0$ & $1$ & $2^s-1$ & $2^{2s}-1$ & $2^{3s}-2^{2s}-2^s+1$\\\hline
 $a\in T_1$& $1$ & $2^s-1$ &  $-1$ & $-2^s+1$\\\hline
 $a\in T_2$ & $1$& $-1$ & $2^s-1$ &   $-2^s+1$\\\hline
 $a\in T_3$& $1$ & $-1$ &$-2^s-1$ & $2^s+1$\\\hline
\end{tabular}
\end{center}
\end{table}
\end{remark}

\subsection{Primitive Association Schemes in $\F_{2^{6s}}$ and $\F_{2^{9s}}$}
In this subsection, we construct  primitive
association schemes in $G=\F_{2^{6s}}$ and $H=\F_{2^{9s}}$. 
\begin{thm}\label{second-ass}
\textup{(i)} Take the following partition of $G$:
\[
R_0'=\{0\},\;R_1'=\bigcup_{i\in T_1} C_i^{(M,G)},\; R_2'=\bigcup_{i\in T_2} C_i^{(M,G)},\;R_3'=\bigcup_{i\in T_3} C_i^{(M,G)}.
\]
Then, $(G,\{R_i'\}_{i=0}^{3})$ is a three-class association scheme, whose parameters are listed in the appendix.

\textup{(ii)} Take the following partition of $H$:
\[
R_0''=\{0\},\;R_1''=\bigcup_{i\in T_1} C_i^{(M,H)},\; R_2''=\bigcup_{i\in T_2} C_i^{(M,H)},\;R_3''=\bigcup_{i\in T_3} C_i^{(M,H)}.
\]
Then, $(H,\{R_i''\}_{i=0}^{3})$ is a three-class association scheme, whose parameters are listed in the appendix.
\end{thm}

\noindent{\bf Proof of Theorem~\ref{second-ass} (i):}
For any $\chi'$ of $G$ such that $\chi'^M=1$, there exists a character $\chi$ of $F^\ast$ such that
\[
\chi|_{E^\ast}=1,\;\chi'=\chi\circ \Norm_{G/F}.
\]
We first compute the Gauss periods
$\eta_a'=\psi'(\gamma^a C_0^{(M,G)})$, $0\le a\le M-1$. By the Hasse-Davenport theorem and
$G_F(\chi)=2^s\sum_{x\in T_1}\chi(\gamma^x)$ (see \cite[Theorem~2.1]{EHKX} or \cite[Lemma~12.0.2]{BEW97} for a proof),   we have
\begin{align*}
\eta_a'&=\frac{1}{M}\sum_{\ell=0}^{M-1}G_G(\chi'^{-\ell})\chi'^{\ell}(\gamma^a)\\
&=-\frac{1}{M}+\frac{-1}{M}\sum_{\ell=1}^{M-1}G_F(\chi^{-\ell})^2\chi^{\ell}(\omega^a)\\
&=-\frac{1}{M}+\frac{-2^s}{M}\sum_{\ell=1}^{M-1}G_F(\chi^{-\ell})\sum_{i \in T_1}\chi^{\ell}(\omega^{a-i})\\
&=-\frac{1}{M}+\frac{-2^s}{M}\left(\sum_{\ell=0}^{M-1}G_F(\chi^{-\ell})\sum_{i \in T_1}\chi^{\ell}(\omega^{a-i})+2^s+1\right)\\
&=-2^s \psi(\omega^a D)-1,
\end{align*}
where $D$ is defined in (\ref{inv-trace0}).
By Lemma~\ref{lem:base}, we obtain
\[
\eta_a'=\begin{cases}2^s-1&\;\textup{ if } a\in T_1,\,\\
-2^{2s}+2^s-1&\;\textup{ if } a\in T_2,\\2^{2s}+2^s-1&\;\textup{ if } a\in T_3.\end{cases}
\]
Now, we compute that
\begin{align*}
\psi'(\ga^a R_k')&=\sum_{i\in T_k}\eta_{a+i}\\
&=(2^s-1)|T_1\cap a+T_k|+(-2^{2s}+2^s-1)|T_2\cap a+T_k|+(2^{2s}+2^s-1)|T_3\cap a+T_k|\\
&=(2^s-1)|T_k|-2^{2s}|T_2\cap a+T_k|+2^{2s}|T_3\cap a+T_k|.
\end{align*}
Clearly, $-2^{2s}|T_2\cap a+T_k|+2^{2s}|T_3\cap a+T_k|$ is the coefficient of $a$ in the element $-2^{2s}(T_2-T_3)T_k^{(-1)}$. The elements $(T_2-T_3)T_k^{(-1)}$, $k=1,2,3$, have been computed in Lemma \ref{gr_t1}.
For each $k=1,2,3$, the sum $\psi(\omega^a R_k')$ follows directly and is listed in Table \ref{tab_2}.
By the Bannai-Muzychuk criterion, $(G,\{R_i'\}_{i=0}^3)$ is a three-class association scheme.
\qed
\vspace{0.3cm}

\begin{table}
\begin{center}
\caption{\label{tab_2}The  values of $\psi(\ga^a R_k')$'s}
\begin{tabular}{|c|c|c|c|c|}\hline
   & $R_0'$& $R_1'$& $R_2'$ &$R_3'$\\\hline
$\ga^a=0$ & $1$ & $(2^{2s}-1)(2^{3s}+1)$ & $2^{s-1}(2^{2s}-1)(2^{3s}+1)$ & $2^{s-1}(2^s-1)^2(2^{3s}+1)$\\\hline
 $a=0$& $1$ & $2^{2s}-1$ &  $2^{s-1}(2^s+1)(-2^{2s}+2^s-1)$ & $2^{s-1}(2^s-1)(2^{2s}+2^s-1)$\\\hline
 $a\in T_1$ & $1$& $-2^{3s}+2^{2s}-1$ & $2^{s-1}(2^{2s}-1)$ &   $2^{s-1}(2^s-1)^2$\\\hline
 $a\not \in T_1\cup \{0\}$& $1$ & $2^{2s}-1$ &$-2^{s-1}$ & $2^{s-1}(-2^{s+1}+1)$\\\hline
\end{tabular}
\end{center}
\end{table}

\begin{remark}
By the proof of Theorem~\ref{second-ass}, the dual scheme of the association scheme in Theorem~\ref{first-ass} is given by
\[
D_0'=\{0\},\;D_1'=C_0^{(M,G)},\; D_2'=\cup_{i\in T_1} C_i^{(M,G)},\;D_3'=\cup_{i\in (T_2\cup T_3) \setminus \{0\}} C_i^{(M,G)}.
\]
Their character values are listed in Table \ref{tab_2d}.
\begin{table}
\begin{center}
\caption{\label{tab_2d}The  values of $\psi(\ga^a D_k')$'s}
\begin{tabular}{|c|c|c|c|c|}\hline
   & $D_0'$& $D_1'$& $D_2'$ &$D_3'$\\\hline
$\ga^a=0$ & $1$ & $(2^{s}-1)(2^{3s}+1)$ & $(2^{2s}-1)(2^{3s}+1)$ & $(2^{2s}-1)^2(2^{2s}-2^s+1)$\\\hline
 $a\in T_1$& $1$ & $2^s-1$ &  $-2^{3s}+2^{2s}-1$ & $(2^s-1)^2(2^s+1)$\\\hline
 $a\in T_2$ & $1$& $-2^{2s}+2^s-1$ & $2^{2s}-1$ &   $-2^s+1$\\\hline
 $a\in T_3$& $1$ & $2^{2s}+2^s-1$ &$2^{2s}-1$ & $-(2^s+1)(2^{s+1}-1)$\\\hline
\end{tabular}
\end{center}
\end{table}
\end{remark}

Next, we give a proof of the second statement of Theorem~\ref{second-ass} (ii).

\noindent{\bf Proof of Theorem~\ref{second-ass} (ii):}
For any multiplicative character $\chi''$ of $H$ such that $\chi''^M=1$, there exists a character $\chi$ of $F^\ast$ such that
\[
\chi|_{E^\ast}=1,\;\chi''=\chi\circ \Norm_{H/F}.
\]
We first compute the Gauss periods
$\eta_{a}''=\psi''(\gamma^a C_0^{(M,H)})$, $0\le a\le M-1$. By the Hasse-Davenport theorem and  $G_F(\chi)=2^s\sum_{i\in T_1}\chi(\omega^i)$, 
we have
\begin{align*}
\eta_a''&=\psi''(\beta^aC_0^{(M,H)})\\
&=\frac{1}{M}\sum_{\ell=0}^{M-1}G_H(\chi''^{-\ell})\chi''^{\ell}(\beta^{a})\\
&=-\frac{1}{M}+\frac{1}{M}\sum_{\ell=1}^{M-1}G_F(\chi^{-\ell})^3\chi^{\ell}(\omega^{a})\\
&=\frac{-1+2^{2s}|T_1|^2}{M}+\frac{2^{2s}}{M}\sum_{\ell=0}^{M-1}G_F(\chi^{-\ell})\sum_{i,j\in T_1}\chi^{\ell}(\omega^{a-i-j}).
\end{align*}
We have $\frac{-1+2^{2s}|T_1|^2}{M}=2^{2s}+2^s-1$. In this case, $\eta_a''$, $0\le a\le M-1$,  take more than three values. Therefore, using Eqn~\eqref{ttt_inv}, we compute directly
\begin{align*}
\psi''(\beta^a R_1'')&-(2^{2s}+2^s-1)|T_1|=\frac{2^{2s}}{M}\sum_{\ell=0}^{M-1}G_F(\chi^{-\ell})\chi^\ell(\omega^a)\sum_{i,j,k\in T_1}\chi^{-\ell}(\omega^{i+j-k})\\
&=\frac{2^{3s}}{M}\sum_{\ell=0}^{M-1}G_F(\chi^{-\ell})\sum_{i \in T_1}\chi^{\ell}(\omega^{-i+a})+\frac{2^{2s}(2^s+1)}{M}\sum_{\ell=0}^{M-1}G_F(\chi^{-\ell})\sum_{i \in \Z_M}\chi^{\ell}(\omega^{i+a})\\
&=2^{3s}\psi(\omega^aD)-2^{2s}(2^s+1).
\end{align*}
Recall that $D=\cup_{i\in -T_1} C_i^{(M,F)}$. It follows that
\[
\psi''(\beta^a R_1'')=2^{2s}-1+2^{3s}\psi(\omega^a D).
\]

By the character values of Lemma~\ref{lem:base}, we obtain
\[
\psi''(\beta^a R_1'')=\begin{cases} -2^{3s}+2^{2s}-1&\;\textup{ if } a\in T_1,\,\\
2^{3s}(2^s-1)+2^{2s}-1&\;\textup{ if } a\in T_2,
\\-(2^s+1)(2^{3s}-2^s+1)&\;\textup{ if } a\in T_3.\end{cases}
\]
In exactly the same way, using Eqn.~\eqref{ttt_inv2} and \eqref{ttt_inv3} we can compute
$\psi''(\beta^a R_k'')$ for $k=2,3$ directly. We record the result in Table \ref{tab_3}. By the Bannai-Muzychuk
criterion, $(H,\{R_i''\}_{i=0}^3)$ is a three-class association scheme, and the above is the first eigenmatrix  of the association scheme; also, the scheme is self-dual.
\qed
\begin{table}
\begin{center}
\caption{\label{tab_3}The  values of $\psi(\beta^a R_k'')$'s}
\begin{tabular}{|c|c|c|c|c|}\hline
   & $R_0''$& $R_1''$& $R_2''$ &$R_3''$\\\hline
$\beta^a=0$ & $1$ & $(2^{2s}-1)(2^{6s}+2^{3s}+1)$ & $2^{s-1}(2^{2s}-1)(2^{6s}+2^{3s}+1)$ & $2^{s-1}(2^s-1)^2(2^{6s}+2^{3s}+1)$\\\hline
 $a\in T_1$& $1$& $-2^{3s}+2^{2s}-1$ & $2^{s-1}(2^{3s}+2^s+1)(2^s-1)$ &  $-2^{s-1}(2^s-1)(2^{3s}-2^s+1)$ \\\hline
 $a\in T_2$ & $1$& $(2^s-1)(2^{3s}+2^s+1)$ & $-2^{s-1}(2^{3s+1}-2^{2s}+1)$ &  $2^{s-1}(2^{s}-1)^2$\\\hline
 $a\in T_3$& $1$ & $-(2^s+1)(2^{3s}-2^s+1)$ &$2^{s-1}(2^{2s}-1)$ & $2^{s-1}(2^{3s+1}+2^{2s}-2^{s+1}+1)$\\\hline
\end{tabular}
\end{center}
\end{table}

\section{Concluding Remarks}
In this note, we gave three constructions of three-class association schemes
from cyclotomy.
In general, one can obtain a three-class association scheme from a two-class association scheme
(a strongly regular Cayley graph) under a certain condition as follows:
Let $X$ be a (additively written) finite abelian group and $S$ be a subset of $X\setminus\{0\}$ such that
$S=-S$. Define $R_0=\{0\}$, $R_1=S$, $R_2=X^\ast\setminus S$.
Assume that $(X,\{R_i\}_{i=0}^{2})$ forms a two-class association scheme, i.e., $\Cay(X,S)$ is strongly regular.
Let $(X,\{D_i\}_{i=0}^2)$ be the dual of $(X,\{R_i\}_{i=0}^2)$, where
we assume that $R_1$ is contained in $D_1$.
Define
\[
R_0'=\{0\}, R_1'=R_1, R_2'=D_1\setminus R_1, R_3'=D_2.
\]
Then, $(X,\{R_i'\}_{i=0}^{3})$ is a three-class association scheme.
This construction is essentially given in \cite[Corollary~3.2]{IM10}.

Our three constructions of three-class association schemes given in this note are not included in this construction. In fact, the association schemes
of Theorems~\ref{first-ass} and \ref{second-ass} (i) are not self-dual but the association scheme obtained from the above construction is self-dual.
Furthermore, neither of the relations of the association scheme
in Theorem~\ref{second-ass} (ii) is strongly regular but two relations of the association scheme from the above construction are strongly regular.

An interesting question that is worth looking into is: what is the relations between the principal part of the first 
eigenmatrices of these three schemes we constructed? According to \cite{bm}, 
the principal parts of the first 
eigenmatrices of the underlying cyclotomic schemes of the same order $M$ satisfy the Hasse-Davenport property, namely, that of $\F_{2^{6s}}$ ($resp.$
$\F_{2^{9s}}$) is square ($resp.$ cube) of that of $\F_{2^{3s}}$ up to a sign. This property seems to be lost after taking fusion using the index sets
$T_1$, $T_2$ and $T_3$.  Our schemes are interesting in that they can serve as examples to test such properties should there
be any.

Finally, we comment that the Gauss periods in consideration takes three values in $\F_{2^{6s}}$ and more in $\F_{2^{9s}}$. It
seems pretty hard to consider the fusions of the cyclotomic scheme in comparison with the case where the Gauss periods take only
two values. Hopefully this will yield us more examples of primitive association schemes with new parameters. We will look into this problem in the future research.
\section*{Appendix: Parameters of the schemes}
Throughout this appendix, we write $q=2^s$. In this appendix, we
give computational results (by Maple) for parameters of the three-class association
schemes we constructed in the previous section. The first and second eigenmatrices of the schemes have been obtained during the proofs.

(1) We have the following computational result for the scheme
$(F,\{R_i\}_{i=0}^3)$ of Theorem~\ref{first-ass}:
Let $A_i$ and $A_i'$ denote the adjacency matrices of $R_i$ and $D_i$ respectively.  With $A_i(A_0,A_1,A_2,A_3)=(A_0,A_1,A_2,A_3)B_i$,  we have
{\tiny
\[
B_1=\begin{pmatrix}
0&q^2-1&0&0\\1&q^2-2&0&0\\0&0&\frac{1}{2} q^2+\frac{1}{2} q-1&\frac{1}{2} q (q-1)\\0&0&\frac{1}{2} q (q+1)&\frac{1}{2} (q-2) (q+1)
\end{pmatrix},
\]
\[
B_2=\begin{pmatrix}
0&0&\frac{1}{2} q^3-\frac{1}{2} q&0\\0&0&\frac{1}{4} q (q^2+q-2)&\frac{1}{4} q^2 (q-1)\\1&\frac{1}{2} q^2+\frac{1}{2} q-1&\frac{1}{4} (q^2+q-6) q&\frac{1}{4} (q^2-3 q+2) q\\0&\frac{1}{2} q (q+1)&\frac{1}{4} (q-2) q (q+1)&\frac{1}{4} (q-2) q (q+1)
\end{pmatrix},
\]
\[
B_3=\begin{pmatrix}
0&0&0&\frac{1}{2} q^3-q^2+\frac{1}{2} q\\0&0&\frac{1}{4} q^2 (q-1)&\frac{1}{4} (q^2-3 q+2) q\\0&\frac{1}{2} q (q-1)&\frac{1}{4} (q^2-3 q+2) q&\frac{1}{4} (q^2-3 q+2) q\\1&\frac{1}{2} (q-2) (q+1)&\frac{1}{4} (q-2) q (q+1)&\frac{1}{4} q (q^2-5 q+6)
\end{pmatrix}.
\]}
With $A_i'(A_0',A_1',A_2',A_3')=(A_0',A_1',A_2',A_3')L_i$, we have
{\tiny
\[
L_1=\begin{pmatrix}
0&q-1&0&0\\1&q-2&0&0\\0&0&0&q-1\\0&0&1&q-2
\end{pmatrix},
\]
\[
L_2=\begin{pmatrix}
0&0&q^2-1&0\\0&0&0&q^2-1\\1&0&q-2&q (q-1)\\0&1&q&-2+q^2-q
\end{pmatrix},
\]
\[
L_3=\begin{pmatrix}
0&0&0&q^3-q^2-q+1\\0&0&q^2-1&q^3-2 q^2-q+2\\0&q-1&q (q-1)&q^3-2 q^2-q+2\\1&q-2&-2+q^2-q&4+q^3-2 q^2-q
\end{pmatrix}.
\]}

(2) We have the following computational result for the scheme
$(G,\{R_i'\}_{i=0}^3)$ of Theorem~\ref{second-ass} (i):
Let $A_i$ and $A_i'$ denote the adjacency matrices of $R_i'$ and $D_i'$ respectively.  With $A_i(A_0,A_1,A_2,A_3)=(A_0,A_1,A_2,A_3)B_i$,  we have
{\tiny
\[
B_1=\begin{pmatrix}
0&q^5+q^2-q^3-1&0&0\\
1&-q^3+q^2+q^4-2&q^3 (q^2-1)/2& (q^2-2 q+1) q^3/2\\
0&(q-1) q^2 (q+1)& (q-1) (q^4+q^3-q^2+q+2)/2&(q-1) q (q^3-q^2-q+1)/2\\
0&(q-1) q^2 (q+1)& q (q^3-q^2-q+1) (q+1)/2& (q^4-2-3 q^3+3 q^2+q) (q+1)/2
\end{pmatrix},
\]
\[
B_2=\begin{pmatrix}
0&0&\frac{1}{2} q^6-\frac{1}{2} q^4+\frac{1}{2} q^3-\frac{1}{2} q&0\\
0&\frac{1}{2} q^3 (q^2-1)&\frac{1}{4} (q^5-2 q^3+2 q^2+q-2) q&\frac{1}{4} (-1+q^4+2 q-2 q^3) q^2\\
1&\frac{1}{2} (q-1) (q^4+q^3-q^2+q+2)&\frac{1}{4} (q^5-3 q^3+3 q^2+q-6) q&\frac{1}{4} q (q^4-2-q^3+3 q) (q-1)\\
0&\frac{1}{2} q (q^3-q^2-q+1) (q+1)&\frac{1}{4} q (q^4-2-q^3+3 q) (q+1)&\frac{1}{4} q (q^4-2-3 q^3+2 q^2+q) (q+1)
\end{pmatrix},
\]
\[
B_3=\begin{pmatrix}
0&0&0&-q^2+\frac{1}{2} q^3-q^5+\frac{1}{2} q+\frac{1}{2} q^4+\frac{1}{2} q^6\\0&\frac{1}{2} (q^2-2 q+1) q^3&\frac{1}{4} (-1+q^4+2 q-2 q^3) q^2&\frac{1}{4} (q^5-4 q^4+6 q^3-2 q^2-3 q+2) q\\0&\frac{1}{2} (q-1) q (q^3-q^2-q+1)&\frac{1}{4} q (q^4-2-q^3+3 q) (q-1)&\frac{1}{4} q (q^4-2-3 q^3+2 q^2+q) (q-1)\\1&\frac{1}{2} (q^4-2-3 q^3+3 q^2+q) (q+1)&\frac{1}{4} q (q^4-2-3 q^3+2 q^2+q) (q+1)&\frac{1}{4} q (q^5+6+7 q^3-q^2-4 q^4-11 q)
\end{pmatrix}.
\]}
With $A_i'(A_0',A_1',A_2',A_3')=(A_0',A_1',A_2',A_3')L_i$, we have
{\tiny
\[
L_1=\begin{pmatrix}
0&q^4-q^3+q-1&0&0\\1&q^2-2&q (q^2-1)&(q^3-2 q^2-q+2) q\\0&q (q-1)&q^2 (q-1)&(q^3-q^2-q+1) (q-1)\\0&(q-2) q&q^3-q^2-q+1&q^4-2 q^3+4 q-2
\end{pmatrix},
\]
\[
L_2=\begin{pmatrix}
 0&0&q^5+q^2-q^3-1&0\\0&q (q^2-1)&q^2 (q^2-1)&q^5-q^4-2 q^3+2 q^2+q-1\\1&q^2 (q-1)&-q^3+q^2+q^4-2&q^2 (q-1) (q^2-1)\\0&q^3-q^2-q+1&q^2 (q^2-1)&q^5-q^4-2 q^3+3 q^2+q-2
\end{pmatrix},
\]
\[
L_3=\begin{pmatrix}
 0&0&0&q^6-q^5+2 q^3-q^2-q^4-q+1\\0&(q^3-2 q^2-q+2) q&q^5-q^4-2 q^3+2 q^2+q-1&q^6-2 q^5-q^4+6 q^3-2 q^2-4 q+2\\0&(q^3-q^2-q+1) (q-1)&q^2 (q-1) (q^2-1)&(q-1) (q^5-q^4-2 q^3+3 q^2+q-2)\\1&q^4-2 q^3+4 q-2&q^5-q^4-2 q^3+3 q^2+q-2&q^6-2 q^5-q^4+6 q^3-4 q^2-6 q+4
\end{pmatrix}.
\]
}

(3)
We have the following computational result  for the scheme
$(H,\{R_i''\}_{i=0}^3)$ of Theorem~\ref{second-ass} (ii):
Let $A_i$ and $A_i'$ denote the adjacency matrices of $R_i''$ and $D_i''$ respectively.
Since this scheme is self dual, we have $Q=P$ (One can check directly that $P^2=q^9I$).
With $A_i(A_0,A_1,A_2,A_3)=(A_0,A_1,A_2,A_3)B_i$,  we have
{\tiny
\[
B_1=\begin{pmatrix}
0&q^8+q^5+q^2-q^6-q^3-1&0&0\\
1&q^7-2 q^5+2 q^4-q^3+q^2-2&\frac{1}{2} (q^5-2 q^3+2 q^2-1) q^3&\frac{1}{2} (q^5-2 q^4+4 q^2-4 q+1) q^3\\
0&(q^5-2 q^3+2 q^2-1) q^2&\frac{1}{2} q^8-q^6+q^5+\frac{1}{2} q^4-\frac{3}{2} q^3+q^2+\frac{1}{2} q-1&\frac{1}{2} (q^7-2 q^6+4 q^4-5 q^3+q^2+2 q-1) q\\
0&(q^5-2 q^3+2 q^2+2 q-1) q^2&\frac{1}{2} (q^7-2 q^5+2 q^4+q^3-3 q^2+1) q&\frac{1}{2} q^8-q^7+2 q^5-\frac{5}{2} q^4-\frac{3}{2} q^3+2 q^2-\frac{1}{2} q-1
\end{pmatrix},
\]
\[
\hspace{-1cm}
B_2=\begin{pmatrix}
0&0&-\frac{1}{2} q-\frac{1}{2} q^4+\frac{1}{2} q^3+\frac{1}{2} q^6-\frac{1}{2} q^7+\frac{1}{2} q^9&0\\0&\frac{1}{2} (q^5-2 q^3+2 q^2-1) q^3&\frac{1}{4} (q^8-2 q^6+2 q^5+q^4-3 q^3+2 q^2+q-2) q&\frac{1}{4} (q^7-2 q^6+4 q^4-5 q^3+q^2+2 q-1) q^2\\1&\frac{1}{2} q^8-q^6+q^5+\frac{1}{2} q^4-\frac{3}{2} q^3+q^2+\frac{1}{2} q-1&\frac{1}{4} (q^8-2 q^6+2 q^5+2 q^4-7 q^3+5 q^2+q-6) q&\frac{1}{4} (q^8-2 q^7+4 q^5-6 q^4+3 q^3+3 q^2-5 q+2) q\\0&\frac{1}{2} (q^7-2 q^5+2 q^4+q^3-3 q^2+1) q&\frac{1}{4} (q^8-2 q^6+2 q^5-3 q^3+3 q^2+q-2) q&\frac{1}{4} (q^8-2 q^7+4 q^5-4 q^4-q^3+5 q^2-q-2) q
\end{pmatrix},
\]
\begin{align*}
&B_3=\\
&\hspace{-1.3cm}\begin{pmatrix}
0&0&0&\frac{1}{2} q-q^2+\frac{1}{2} q^4+\frac{1}{2} q^3-q^5+\frac{1}{2} q^6-q^8+\frac{1}{2} q^7+\frac{1}{2} q^9\\0&\frac{1}{2} (q^5-2 q^4+4 q^2-4 q+1) q^3&\frac{1}{4} (q^7-2 q^6+4 q^4-5 q^3+q^2+2 q-1) q^2&\frac{1}{4} (q^8-4 q^7+6 q^6-2 q^5-7 q^4+9 q^3-2 q^2-3 q+2) q\\0&\frac{1}{2} (q^7-2 q^6+4 q^4-5 q^3+q^2+2 q-1) q&\frac{1}{4} (q^8-2 q^7+4 q^5-6 q^4+3 q^3+3 q^2-5 q+2) q&\frac{1}{4} (q^8-4 q^7+6 q^6-2 q^5-6 q^4+9 q^3-3 q^2-3 q+2) q\\1&\frac{1}{2} q^8-q^7+2 q^5-\frac{5}{2} q^4-\frac{3}{2} q^3+2 q^2-\frac{1}{2} q-1&\frac{1}{4} (q^8-2 q^7+4 q^5-4 q^4-q^3+5 q^2-q-2) q&\frac{1}{4} (q^8-4 q^7+6 q^6-2 q^5-8 q^4+13 q^3+3 q^2-11 q+6) q
\end{pmatrix}.
\end{align*}}
In this case, we have $A_i=A_i'$ and $B_i=L_i$ for $0\le i\le 3$.


\end{document}